\documentclass{article}

\usepackage{amsthm}

\newtheorem{conj}{Conjecture}
\newtheorem{lemma}{Lemma}
\newtheorem{theorem}{Theorem}
\newtheorem{corollary}{Corollary}

\theoremstyle{definition}
\newtheorem{remark}{Remark}
\newtheorem{example}{Example}

\usepackage{amssymb}
\usepackage{amsmath}

\newcommand{\N}{\mathbb{N}}
\newcommand{\R}{\mathbb{R}}
\newcommand{\PP}{\mathbb{P}}
\newcommand{\F}{{\cal F}}
\newcommand{\E}{{\bf E}}
\renewcommand{\a}{\alpha}
\renewcommand{\i}{\infty}
\newcommand{\e}{\varepsilon}
\def\as{a.s.~}

\newcommand{\keywords}{{\em Keywords:}}

\begin{document}


\title{Almost Sure Convergence of Solutions to Non-Homogeneous
  Stochastic Difference Equation}
\author{Gregory Berkolaiko$\dagger$ and Alexandra Rodkina$^\ddagger$\\
  $^{\dagger}$Department of Mathematics, Texas A\&M University,\\ College
  Station, U.S.A.\\
  $^{\ddagger}$Department of Mathematics and Computer Science,\\
  University of the West Indies, Kingston, Jamaica}


\markboth{Non-homogeneous stochastic difference equations}{G.~Berkolaiko and
  A.~Rodkina}

\maketitle

\begin{abstract}
  We consider a non-homogeneous nonlinear stochastic difference
  equation
  \begin{displaymath}
    X_{n+1} = X_n \Big(1+f(X_n)\xi_{n+1}\Big) + S_n,
    \quad n=0,1,\ldots,
  \end{displaymath}
  and its linear counterpart
  \begin{displaymath}
    X_{n+1} = X_n \Big(1+\xi_{n+1}\Big) + S_n,
    \quad n=0,1,\ldots,
  \end{displaymath}
  both with initial value $X_0$, non-random decaying free
  coefficient $S_n$ and independent random variables $\xi_n$. We
  establish results on \as convergence of solutions $X_n$ to zero.
  Obtained necessary conditions tie together certain moments of
  the noise $\xi_n$ and the rate of decay of $S_n$. To ascertain
  sharpness of our conditions we discuss some situations when
  $X_n$ diverges. We also establish a result concerning the rate
  of decay of $X_n$ to zero.

  Several examples are given to illustrate the ideas of the paper.
\end{abstract}

\hspace{1cm}

{\keywords Nonlinear stochastic difference equations, almost
 sure stability, martingale convergence theorem.

}

\section{Introduction}

The theory of stochastic difference equations is relatively young,
especially in its nonlinear part. Linear stochastic difference
equations with independent identically distributed perturbations
(i.i.d.) are the most studied ones (cf~\cite{K}) but even for this
type of equation there still exist some open questions \cite{RB}.
In this paper we are going to give answers to some of them and
then proceed to discuss a class of nonlinear stochastic difference
equations for which very few results are available \cite{KS, KS02,
Rodkina01, RMK, RS, RS1}.

The interest towards stochastic difference equations has been on the
increase due to their numerous applications and the fact that they
serve for numerical simulations of stochastic differential
equations (cf~\cite{Higham, HMS, KP, Mao97}). Stability of solutions of
stochastic difference equations is also very important in, to give
some examples, mathematical finance (asset price evolution in
discrete $(B,S)$-markets) and mathematical biology (population
dynamics), see, for example, \cite{FH} and references therein.

The main objects of our consideration are the following equations:
the non-homogeneous nonlinear stochastic difference equation
\begin{equation}
  \label{eq:nlin_recng}
  X_{n+1} = X_n \Big(1+f(X_n)\xi_{n+1}\Big) + S_n,\qquad 
  n\in\N_0 = \N\cup\{0\},
\end{equation}
and its linear counterpart
\begin{equation}
  \label{eq:lin_recng}
  X_{n+1} = X_n \Big(1+\xi_{n+1}\Big) + S_n, \qquad n\in\N_0,
\end{equation}
with initial value $X_0>0$, non-random free coefficient $S_n$ and
independent random variables $\xi_n$.  Unless explicitly indicated, we
do not demand that $\xi_n$ be identically distributed.  Everywhere in
the paper we suppose that
\begin{eqnarray}
  \label{eq:f_cond}
  &&f:\R^1\to [0, 1] \mbox{ is continuous and }
  f(u)=0 \Leftrightarrow u=0,\\
  \nonumber
  &&1+\xi_{n+1}>0 \quad \mbox{and} \quad S_n>0
  \qquad \forall \, n\in\N_0.
\end{eqnarray}
These conditions guarantee that $X_n$ remains positive for all
$n$.

Equations of the type (\ref{eq:nlin_recng}) and (\ref{eq:lin_recng}) are
sufficiently complex to require more powerful methods than those used to
study, for example, the linear homogenous equation
\begin{equation}
  \label{eq:lin_homo}
  X_{n+1} = X_n \Big(1+\xi_{n+1}\Big), \qquad n\in\N_0.
\end{equation}
On the other hand, equations (\ref{eq:nlin_recng}) and (\ref{eq:lin_recng})
are sufficiently simple to allow a rather complete understanding of their
behaviour.  In our paper we use an adaptation of a martingale convergence
theorem to prove most of the results.  The methods of proof that we develop
can also be used on more complicated recursions or in more applied contexts,
for example to study the faithfulness properties\footnote{such as the
  A-stability, which was studied in \cite{Higham} on the example of the
  recursion of the type (\ref{eq:lin_homo}).  A method is said to be {\em
    A-stable} if it correctly predicts the asymptotic stability of the
  approximated equation.}  of numerical solutions to stochastic differential
equations.

To get the flavour of our results it is instructive to start with
the behaviour of the corresponding deterministic equation,
\begin{displaymath}
  x_{n+1} = x_n \Big(1+a_{n+1}\Big) + S_n,
\end{displaymath}
with $1+a_{n+1}>0$ (the nonlinear deterministic equation is discussed in
Section~\ref{sec:det_lemma}). If $a_n\equiv a$, the solutions converge to zero
when $a<0$ (or $\ln(1+a)<0$, which is the same) and $S_n\to 0$.

Now take $S_n\equiv0$ and allow $a_n$ to contain noise,
$a_n=a+\zeta_n$ with $\E\zeta_n=0$.  It is easy to see that the
solutions will still tend to zero if $a = \E(a+\zeta_n) < 0$.  But
they will also tend to zero if $a>0$ but $\E\ln(1+a+\zeta_n)<0$, which
is a weaker condition.  We will refer to this phenomenon as the
``stabilisation by noise'': the solution of $x_{n+1} = x_n (1+a)$ with
$a>0$ can be stabilised by adding some noise to $a$ (for an in-depth
discussion of stabilisation by noise see
e.g.~\cite{JAXM:05a,Arn,Car:2,M94b}).  A natural question arises: when
the noise is present, how fast must $S_n$ decay to guarantee that the
convergence persists?  Would $S_n\to0$ be enough?  We will discuss
this question at length in the present paper but the short answer is
the following.  The coefficients $S_n$ must have a power law decay,
with the exponent determined by the nature of the noise.  Thus, the
addition of the noise stabilises the homogenous linear equation but
imposes stronger conditions on the free coefficient $S_n$ of the
non-homogenous one.  It is interesting to compare our results with
those available in the continuous case, where the interplay of the
noise and the rate of decay of the free coefficient was studied in
\cite{AR}.

In nonlinear case (\ref{eq:nlin_recng}), however, the noise does
not have such stabilising effect. Our stability result (if
restricted to i.i.d.~noises) includes only the case $\E\xi_n<0$.
We investigate the case $\E\xi_n>0$ further and show, for bounded
i.i.d.~$\xi_n$, that $\lim_{n\to \infty}X_n=0$ with probability
zero. Heuristically, the noise does not have the stabilising
effect on the nonlinear equation because the coefficient by
$\xi_n$ becomes too small if $X_n\to0$ (see
condition~(\ref{eq:f_cond})). The situation changes when instead
of equation (\ref{eq:nlin_recng}) we consider a discrete version
of Ito stochastic equation with the drift and diffusion parts
separated and multiplied by coefficients with different scaling:
\begin{equation}
  \label{eq:discr_Ito}
  X_{n+1} = (1 + kf(X_n)a+\sqrt{kf(X_n)}\zeta_{n+1}) X_n+S_n,
  \quad n\in\N_0.
\end{equation}
In this case we give a sufficient conditions for $\lim_{n\to
\infty}X_n=0$ \as even when $a$ is positive (but not too large).

The structure of the paper is as follows. In section \ref{sec:aux}
we give some necessary definitions and state two lemmas.
Lemma~\ref{lem:nonegdif} can be considered a discrete version of
martingale convergence theorem and is the main tool we use to
prove our results. Section \ref{sec:linind} is devoted to the
\as convergence to zero of solutions to the linear
equations with independent noises. We also present a result on the
rate of decay of the solutions. Section \ref{sec:disc_lin} is
devoted to a discussion of the obtained results as they apply to
the i.i.d.~noise. Further results in this simple case highlight
some aspects of behaviour of the solutions. In particular, we
construct some examples that indicate that our conditions for the
\as convergence might be necessary as well as sufficient. We also
find that when the decay of $S_n$ is insufficient to guarantee
convergence but $\E\ln(1+\xi_n)$ is negative, the lower limit of
the solution is still zero. This implies that, in some cases, the
solution will oscillate with increasing amplitude.

Section \ref{sec:nonlin} is devoted to nonlinear equation
(\ref{eq:nlin_recng}). Sufficient conditions which guarantee that $\lim_{n\to
  \infty}X_n=0$ are given in the case when $S_n$ are summable and when
$S_n^\a$ are summable with some $\a<1$. We also prove, for bounded
i.i.d.~$\xi_n$ with $\E\xi_n>0$, that $\lim_{n\to \infty}X_n=0$ with
probability zero. Then we consider equation (\ref{eq:discr_Ito}), a discrete
version of Ito stochastic equation, and give a sufficient conditions for the
\as convergence of the solutions to zero.

We illustrate our results with examples and defer all
proofs to the last section of the paper.

\section{Auxiliary Definitions and Facts}\label{sec:aux}

Let $(\Omega, {\F}, \{{\F}_n\}_{n \in \N}, {\PP})$ be a complete
filtered probability space.  Let $\{\xi_{i}\}_{i\in\N}$ be a
sequence of independent random variables. We suppose that
filtration $\{{\F}_n\}_{n \in \N}$ is naturally generated:
$\F_{n+1} = \sigma \{\xi_{i+1}: i\leq n\}$.  Among all sequences
$\{X_n\}_{n\in\N}$ of random variables we distinguish those for
which $X_n$ are ${\F}_n$-measurable $\forall n \in \N$.

We use the standard abbreviation ``{\it a.s.}'' for the wordings
``almost sure'' or ``almost surely'' with respect to the fixed
probability measure ${\PP}$ throughout the text.

A stochastic sequence $\{X_n\}_{n \in N}$ is said to be an {\em
  $\F_n$-martingale}, if ${\bf E}|X_n|<\i$ and ${\bf E}\bigl(X_n\bigl|
{\cal F}_{n-1}\bigr)=X_{n-1}$ \as for all $n\in\N$. A stochastic
sequence $\{\mu_n\}_{n\in\N}$ is said to be an {\em
  $\F_n$-martingale-difference}, if ${\bf E}|\mu_n|<\i$ and ${\bf
    E}\bigl(\mu_n\bigl| \F_{n-1}\bigr)=0$ \as for all
  $n\in\N$.

For more details on stochastic concepts and notation we refer the
reader to \cite{LSh, Mao97, Nev, Shiryaev96}.

Below is a version of a martingale convergence theorem, which is
convenient for many proofs.

\begin{lemma}\label{lem:nonegdif}
  Let $\{Z_n\}_{n\in \N}$ be a non-negative  $\F_n$-measurable
  process, ${\E}|Z_n|<\i$  $\forall n\in \N$ and
  \begin{displaymath}
    Z_{n+1}\le Z_{n}+u_n-v_n+\nu _{n+1}, \quad n\in\N,
  \end{displaymath}
  where $\{\nu_n\}_{n\in \N}$ is $\F_n$-martingale-difference,
  $\{u_n\}_{n\in \N}$, $\{v_n\}_{n\in \N}$ are nonnegative ${\cal
    F}_{n}$-measurable processes, ${\E}|u_n|$ and ${\E}|v_n|$ are finite.

  Then
  $$
  \left\{\omega: \sum_{n=1}^{\i} u_n<\i\right\}
  \subseteq \left\{\omega: \sum_{n=1}^{\i} v_n<\i\right\}\bigcap\{Z\to\}.
  $$
\end{lemma}
Here $\left\{Z\to \right\}$ denotes the set of all
$\omega\in\Omega$ for which $Z_\infty=\lim\limits_{n\to\infty}Z_n$
exists and is finite.

We will also use the following elementary estimate.
\begin{lemma}
  \label{lem:power_ineq}
  For any $\alpha\geq1$ there exists a function $K$ continuous
  on $(0,\infty)$ such that for any $a>0$ and $b>0$
  \begin{displaymath}
    (a+b)^\alpha \leq (1+\epsilon)a^\alpha + K(\epsilon)b^\alpha,
  \end{displaymath}
  where $K(\e)$ can be estimated in the following way:
  \begin{displaymath}
    K(\e) \le 1 + K_1(\alpha)\e^{1-\alpha}.
  \end{displaymath}
\end{lemma}
We define $[u]^+$ and $[u]^-$ to be the positive and negative parts of $u$
correspondingly,
\begin{displaymath}
  [u]^+= \left\{
    \begin{array}{ll}
      u,&{\rm if}\:u>0,\\ 0, &{\rm otherwise,}
    \end{array}\right. \qquad
  [u]^-= \left\{
    \begin{array}{ll}
      u,&{\rm if}\:u<0,\\ 0, &{\rm otherwise.}
    \end{array}\right.
\end{displaymath}
We will say that a sequence $\{S_n\}$ is $\alpha$-summable if
\begin{displaymath}
  \sum_{n=1}^\infty S_n^\alpha < \infty.
\end{displaymath}

\section{Linear non-homogeneous equation with independent noises.}
\label{sec:linind}

Below is our main result on the limit of solutions to linear
equation (\ref{eq:lin_recng}).  The conditions for \as existence of
a limit depend on the balance between $\alpha$-summability of $S_n$ and the
signs of $\E(1+\xi_{i+1})^\alpha-1$.
\begin{theorem}
  \label{thm:indlin}
  Let $X_n$ be a solution to equation (\ref{eq:lin_recng}).
  If there exists $\alpha>0$ such that
  \begin{equation}
    \label{eq:pos_part_conv}
    \sum_{i=1}^{\infty}\left[\E(1+\xi_{i+1})^\alpha-1\right]^+<\infty,
  \end{equation}
  and
  \begin{eqnarray}
    \label{eq:S_alpha_sum}
    &&\sum_{i=1}^{\infty}S_i^\alpha<\infty,
    \qquad \qquad \mbox{if}\quad \alpha\leq1,\\
    \label{eq:S_alpha_sum_plus}
    &&\sum_{i=1}^{\infty}\frac{S_i^\alpha}
    {\bigl|1-\E(1+\xi_{i+1})^\alpha\bigr|^{\alpha-1}} < \infty,
    \qquad\qquad \mbox{if}\quad\alpha>1,
  \end{eqnarray}
  then $\lim_{n\to\infty} X_{n}$ exists. If, in addition,
  \begin{equation}
    \label{eq:neg_part_div}
    \sum_{i=1}^{\infty}[\E(1+\xi_{i+1})^\alpha-1]^-=-\infty,
  \end{equation}
  then $\lim_{n\to\infty} X_{n}=0$.
\end{theorem}

\begin{remark}
  If $\xi_n$ are independent and identically distributed (i.i.d.),
  as opposed to just independent, then $\E(1+\xi_{n+1})^{\a}-1$
  does not depend on $n$. Therefore conditions
  (\ref{eq:pos_part_conv}) and (\ref{eq:neg_part_div})
  are fulfilled whenever $\E(1+\xi_{n+1})^{\a}-1 < 0$ for the
  corresponding value of $\alpha$.

  We note that if $\beta<\alpha$ then
  $\E(1+\xi_{n+1})^\a - 1 < 0$ implies
  $\E(1+\xi_{n+1})^\beta - 1 < 0$.  Thus the requirements
  on $\xi$ get stronger with the growth of $\alpha$.
  This is compensated by weakening of the requirements on $S_n$  (in
  the i.i.d. case condition (\ref{eq:S_alpha_sum_plus}) is just the
  $\alpha$-summability of $S_n$).

  Interestingly, when $\alpha<1$ one can have $\E\xi_n>0$.
  This will be discussed in more detail in Section \ref{sec:comp}
  below.
\end{remark}

The following example illustrates the case when
$\sum_{i=1}^{\i}S_i=\i$ and $\alpha>1$.

\begin{example}
  Let
  \[
  \xi_n=\left\{\begin{array}{cc}
      -n^{-\frac 13} \quad &\mbox{with probability} \quad 1-\frac 1{n^2},\\
      \sqrt{n} \quad &\mbox{with probability} \quad \frac 1{n^2},
    \end{array}\right.
  \]
  and
  \begin{displaymath}
    S_n\sim n^{-\frac34}.
  \end{displaymath}

  Then
  \[
  \E\xi_n = n^{-\frac 13}\left(1-\frac 1{n^2}\right) 
  + \sqrt{n}\frac1{n^2}\sim -n^{-\frac 13},
  \]
  and
  \[
  \E\xi_n^2=-n^{-\frac 23}\left(1-\frac 1{n^2}\right)+n\frac
  1{n^2}\sim n^{-\frac 23}.
  \]
  Therefore,
  \[
  1-\E(1+\xi_n)^2=-2\E\xi_n-\E\xi_n^2\sim 2n^{-\frac 13}.
  \]
  Even though $S_n$ are not summable, conditions
  (\ref{eq:neg_part_div}) and (\ref{eq:S_alpha_sum_plus})
  are fulfilled with $\a=2$, since
  \[
  \sum_{n=1}^{\i}[\E(1+\xi_{n+1})^2-1]\sim
  -2\sum_{n=1}^{\i}n^{-\frac 13}=-\i,
  \]
  \[
  \sum_{n=1}^{\i}\frac{S_n^{2}}{1-\E(1+\xi_{n+1})^{2}}\sim
  \sum_{n=1}^{\i}n^{\frac 13}n^{-\frac 64}=\sum_{n=1}^{\i}n^{-\frac
    76}<\i.
  \]
  Then Theorem \ref{thm:indlin} implies that
  $\lim_{n\to\infty} X_{n}=0$ \as
\end{example}

The next result gives the rate of decay of solutions to
equation (\ref{eq:lin_recng}) when we impose more restriction on the
summability of the free coefficient $S_n$.
\begin{theorem}
  \label{thm:indexp}
  Let $\xi_n$ be independent random variables and
  $X_n$ be a solution to equation (\ref{eq:lin_recng}).  If for some
  $\a \in (0,1]$ there are $\kappa_i$ such that
  \begin{equation}
    \label{eq:kappa1}
    \kappa_i \ge  \bigl[\E(1+\xi_{i+1})^\a-1\bigr]^-,
  \end{equation}
  \begin{equation}
    \label{eq:kappa2}
    \sum_{i=1}^{\i}\kappa_i = -\i,
  \end{equation}
  \begin{displaymath}
    \sum_{n=1}^{\i}e^{-\sum_{i=1}^{n+1}\kappa_i}S_n^\a<\i,
  \end{displaymath}
  then for every $\gamma\in (0, 1)$
  \begin{displaymath}
    \lim_{n\to\infty} e^{-\gamma\sum_{i=1}^n\kappa_i}X_n^\a=0.
  \end{displaymath}
\end{theorem}

\section{Discussion of Theorem \ref{thm:indlin}}
\label{sec:disc_lin}

In this section we limit ourselves to considering i.i.d.~$\xi_n$.
We discuss two questions here, the sharpness of the conditions of
Theorem~\ref{thm:indlin} and using $\E\ln\big(1+\xi_{i}\big)<0$ as
an indicator of \as convergence.

\subsection{Is $\a$-summability necessary?}

The following lemma shows that in general one can not relax the condition of
$\alpha$-summability of $S_n$.
\begin{lemma}\label{lem:limsup_counter}
  For any $\alpha$ and $\beta$ satisfying $0<\alpha<\beta$ there exist
  i.i.d.~random variables $\{\xi_n\}_{n=1}^\infty$ and perturbations $S_n$
  such that
  \begin{eqnarray}
    \label{eq:limsup_Ealpha}
    &&\E(1+\xi)^\alpha = 1,\\
    \label{eq:limsup_betaSum}
    &&\sum_{n=1}^\infty S_n^\beta < \infty,
  \end{eqnarray}
  and yet the solution $X_n$ of equation~(\ref{eq:lin_recng}) is diverging in
  the sense that
  \begin{displaymath}
    \limsup_{n\to\infty} X_n = \infty \qquad \mbox{\as}
  \end{displaymath}
\end{lemma}

\begin{remark}\label{rem:limsup_counter}
  Theorem \ref{thm:indlin} requires $\alpha > \beta$ to guarantee
  \as convergence of $X_n$.
\end{remark}

\subsection{Homogeneous equation}

When equation (\ref{eq:lin_recng}) is homogeneous ({\it i.e.}
$S_n=0$), the limit is zero if and only if
$\E\ln\big(1+\xi_{i}\big)<0$ (see, for example, \cite{Nev} or
\cite{RB}):
\begin{theorem}
  \label{thm:iidhom}
  Assume that $\{\xi_n\}_{n\in\N}$ are
  i.i.d.~random variables and $X_n$ is the solution of equation
  (\ref{eq:lin_recng}) with $S_n=0$.  Then
  $\lim_{n\to+\i} X_n = 0$ \as if and only if
  $\E\ln\big(1+\xi_{i}\big)<0$.
\end{theorem}
It seems, therefore, that $\E\ln\big(1+\xi_{i}\big)<0$ is a
natural indicator of the convergence of $X_n$.  Sections
\ref{sec:liminf} and \ref{sec:comp} develop this observation.

\subsection{Lower limit}
\label{sec:liminf}

When $\E\ln\big(1+\xi_{i}\big)<0$ and $S_n$ is non-zero but decreases
exponentially with $n$, it was proved in \cite{RB} that $\lim_{n \to +\i} X_n
= 0$.  When $S_n$ does not decrease as rapidly, it turns out that condition
$\E\ln\big(1+\xi_{i}\big)<0$ guarantees that the lower limit of $X_n$ is equal
to zero.
\begin{theorem}
  \label{thm:liminf}
  Let $\xi_n$ be i.i.d. with $\E\ln(1+\xi_{n+1}) < 0$.  If there
  exists $\alpha>0$ such that $\sum_{i=1}^{\i}S_i^{\a}<\i$, then
  \begin{displaymath}
    \liminf_{n\to\infty} X_{n} = 0 \qquad \as
  \end{displaymath}
\end{theorem}

\begin{remark}
  In some cases, in particular those covered by
  Lemma~\ref{lem:limsup_counter}, the lower limit is equal to zero while the
  limit does not exist.  An interesting question is the existence of the
  limiting distribution of $X_n$ in such cases.
\end{remark}

\subsection{Connection between $\E\ln(1+\xi_{i})$ and $\E(1+\xi_i)^{\a}-1$.}
\label{sec:comp}

Theorem \ref{thm:iidhom} indicates that the sign of
$\E\ln(1+\xi_{i})$ is crucial in the question of stability of
homogenous equation with i.i.d. noises.  Theorem~\ref{thm:indlin}, however,
depends on the sign of $\E(1+\xi_i)^\alpha-1$ to establish stability.  The
following lemma provides the connection between the two expectations.
\begin{lemma}\label{lem:alnind0}
  Let $\xi$ be such that $\PP(\xi>0) > 0$.  Then $\E\ln(1+\xi)<0$
  if and only if there exists $\alpha>0$ such that $\E(1+\xi)^\alpha-1=0$.
  If such $\alpha$ exists then
  \begin{equation}
    \label{eq:all_beta_less0}
    \E(1+\xi)^\beta-1<0 \qquad \forall\beta\in(0,\alpha).
  \end{equation}
\end{lemma}
\begin{proof}
  The harder ``only if'' part was proved in \cite{K}, using that
  $\E(1+\xi)^\alpha$ is a convex function and its derivative at $\alpha=0$
  is equal to $\E\ln(1+\xi)$.  Convexity also implies
  inequality (\ref{eq:all_beta_less0}).

  To prove the ``if'' part we take expectation of the both parts of the
  inequality $(1+\xi)^u \geq 1+ u\ln(1+\xi)$ which can be obtained by
  truncating the Taylor series of $(1+\xi)^u$ with respect to $u$.
\end{proof}

When $\xi_n$ are not identically distributed, one needs a uniform bound on
$\alpha$.  Such a bound is provided by the following lemma.

\begin{lemma}\label{lem:alnind}
  Suppose that there exists some constant $K>0$ such that
  \begin{equation}
    \label{eq:ln2ln}
    \frac{\E \left[(2+\xi_i) \ln^2(1+\xi_i)\right]}
    {\left|\E \ln(1+\xi_i)\right|} \le K, \quad \forall\, i\in \N.
  \end{equation}
  Then for all $\a$ satisfying
  \begin{displaymath}
    \a < \min\left(1/K, 1\right)
  \end{displaymath}
  one has
  \begin{equation}
    \label{eq:uniform_bound}
    \a \E\ln(1+\xi_i) \le \E(1+\xi_i)^\alpha - 1
    \le \alpha \left(\E\ln(1+\xi_i) + \frac{|\E\ln(1+\xi_i)|}{2}\right).
  \end{equation}
\end{lemma}

\begin{example}
  Suppose that $-1<-k\le \xi_n\le L$
  and $\left|\E\ln(1+\xi_n)\right|\ge c$ for some $k,L,c>0$ uniformly
  in $n\in\N$. Then condition (\ref{eq:ln2ln}) is fulfilled.
\end{example}

\subsection{Reformulation of Theorem~\ref{thm:indlin} in the i.i.d. case}

Following the discussion of the previous sections we can reformulate
Theorem~\ref{thm:indlin} in this concise way.

\begin{corollary}
  Let $\xi_n$ be i.i.d.~random variables satisfying
  \begin{displaymath}
    \E(1+\xi_n)^\alpha-1 \leq 0
  \end{displaymath}
  for some $\alpha>0$.
  If $S_n$ are $\alpha$-summable then the solutions of
  \begin{displaymath}
    X_{n+1} = X_n \Big(1+\xi_{n+1}\Big) + S_n, \qquad n\in\N_0,
  \end{displaymath}
  converge to zero \as
\end{corollary}

\begin{proof}
  The only part of the statement that does not obviously follow
  from Theorem~\ref{thm:indlin} is what happens when
  $\E(1+\xi_n)^\alpha-1 = 0$. In this case
  Theorem~\ref{thm:indlin} guarantees only the existence of a
  limit. Here, however, we employ Lemma~\ref{lem:alnind0} to infer
  that $\E\ln(1+\xi_n) < 0$. Then we use Theorem~\ref{thm:liminf}
  to confirm that the limit must indeed be zero.
\end{proof}

On the other hand, Lemma~\ref{lem:limsup_counter} hints that
$\alpha$-summability is not only a sufficient but also a necessary
condition. We formulate this guess as a conjecture.
\begin{conj}
  Let $\xi_n$ be i.i.d.~random variables satisfying
  $\E\ln(1+\xi_n) < 0$ and let $\alpha>0$ be such that
  \begin{displaymath}
    \E(1+\xi_n)^\alpha - 1 = 0.
  \end{displaymath}
  Then the solutions of
  \begin{displaymath}
    X_{n+1} = X_n \Big(1+\xi_{n+1}\Big) + S_n, \qquad n\in\N_0
  \end{displaymath}
  \as converge to zero {\bf if and only if} $S_n$ is
  $\alpha$-summable.
\end{conj}

Another interesting question would be to study the convergence of
$X_n$ to zero {\em in probability}.  For this type convergence, it
might be possible to relax the conditions on the decay of $S_n$.
Previous results by various authors \cite{KesSpi84,Bou87,Muk91} should
be helpful in this direction.

\section{Nonlinear equation}
\label{sec:nonlin}

In this section we consider nonlinear recursion of the type
\begin{equation}
  \label{eq:nonlin_eq}
  X_{n+1} = X_n \Big(1+f(X_n)\xi_{n+1}\Big) + S_n,
  \qquad n\in\N_0
\end{equation}
with independent random variables $\xi_n$. As mentioned earlier we
assume that the function $f(u)$ is continuous with values in the
interval $[0,1]$ and is equal to 0 only at $u=0$. We also assume
that both terms in equation (\ref{eq:nonlin_eq}) are non-negative for all
$n$.

\subsection{Convergence results for nonlinear equation with independent
  noises}
\label{sec:nonlin_thms}

Only the $\alpha=0$ case of Theorem~\ref{thm:indlin} really
carries over to the nonlinear equations of
type~(\ref{eq:nonlin_eq}).
\begin{theorem} \label{thm:posnegn}
  If the components of equation (\ref{eq:nonlin_eq}) satisfy the conditions
  detailed above, $S_n$ are summable and
  \begin{equation}
    \label{eq:S_n_and_pos_exp_conv}
    \sum_{n=1}^\infty \left[\E\xi_n\right]^+ < \infty,
  \end{equation}
  then $\lim_{n\to\infty} X_{n}$ exists.  If, in addition,
  \begin{equation}
    \label{eq:neg_exp_div}
    \sum_{n=1}^\infty \left[\E\xi_n\right]^- = -\infty,
  \end{equation}
  then $\lim_{n\to\infty} X_{n}=0$.
\end{theorem}
In the case when $S_n$ are $\alpha$-summable for some $\a\in (0,
1)$ we obtain a much more restrictive result compared with Theorem
\ref{thm:indlin}.  Even the case of i.i.d.~noises with positive
$\E\xi_n$ is not covered by this theorem.  We will explore the
reason for this in the next section.
\begin{theorem}
  \label{thm:asumind}
  Let $S_n$ be $\alpha$-summable for some $\a\in(0,1)$,
  \begin{equation}
    \sum_{n=1}^{\i}[\E\xi_n]^+<\i,
    \label{eq:conv13}
  \end{equation}
  then $\lim_{n\to\infty} X_{n}$ exists.  If, in addition, $3\E\xi_n^2 -
  (2-\a)[\E\xi_n^3]^+ > 0$ starting with some $n$ and
  \begin{equation}\label{eq:conv2}
    \sum_{n=1}^{\i}\bigg(\E\xi_n^2-\frac{2-\a}{3}[\E\xi_n^3]^+\bigg)
    = \i,
  \end{equation}
  then $\lim_{n\to\infty} X_{n}=0$.
\end{theorem}
The following example shows that in the case when $S_n^\alpha$ are
summable with some $\a<1$, Theorem \ref{thm:asumind} gives less
restrictive conditions than Theorem \ref{thm:posnegn}.
\begin{example}
  Let $\xi_n$ be uniformly distributed on the interval $[-1+n^{-2}, 1]$.
  Then
  \begin{displaymath}
    \E\xi_n \propto n^{-2},\qquad
    \E\xi_n^2 \propto 1 \qquad \mbox{ and }
    \qquad \E\xi_n^3 \propto n^{-2}.
  \end{displaymath}
  Thus conditions~(\ref{eq:conv13}) and~(\ref{eq:conv2}) are fulfilled
  for all $\a\in (0,1)$, but condition~(\ref{eq:neg_exp_div}) is not.
\end{example}

\subsection{Divergence in nonlinear equation with $\E\xi_i>0$}
\label{sec:nonlin_pos_mean}

In this subsection we present a result explaining why one cannot
fully generalise Theorem~\ref{thm:indlin} to nonlinear equations
of the type~(\ref{eq:nonlin_eq}).

\begin{theorem} \label{thm:div}
  Let $X_n$ be a solution of equation (\ref{eq:nonlin_eq}) with i.i.d.~$\xi_n$
  satisfying
  \begin{displaymath}
    \E\xi_n>0 \qquad\mbox{and}\qquad -1<-k_0 \le \xi_n \le L, \qquad n\in \N.
  \end{displaymath}
  Then $\PP\{X_n\to 0\} = 0$.
\end{theorem}

However, if, instead of equation (\ref{eq:nonlin_eq}) we consider a discrete
analogue of Ito equation, the situation is reversed and we obtain a
convergence result when $S_n^\a$ is summable with some $\a>0$.

\subsection{Analogue of Ito equation}
\label{sec:nonlin_ito}

We consider the  discrete analogue of Ito equation
\begin{equation}
  \label{eq:nlinIto}
  X_{n+1} = (1 + kf(X_n)a+\sqrt{kf(X_n)}\zeta_{n+1}) X_n+S_n, \quad X_0>0,
\end{equation}
where $a>0$, $\zeta_n$ are i.i.d., $\E\zeta_n=0$, $\E\zeta_n^2<\i$
and $\E|\zeta_n|^3<\i$.

We assume, as everywhere before, that for all positive $u$ and all $n$
\begin{equation}
  \label{eq:Ito_pos}
  1 + kf(u)a+\sqrt{kf(u)}\zeta_{n+1} > 0
  \qquad \mbox{and} \qquad
  S_n \ge 0.
\end{equation}
\begin{theorem}
  \label{thm:Ito}
  Let conditions (\ref{eq:f_cond}) and (\ref{eq:Ito_pos}) be fulfilled. Let
  also
  \begin{displaymath}
    a < \frac{\E \zeta^2}{2}.
  \end{displaymath}
  Suppose $S_n$ are $\alpha$-summable for some $\a$ satisfying  the
  inequality
  \begin{equation}
    \label{eq:Ito_alpha}
    \a<\a_0=\frac{\E \zeta^2-2a}{\E \zeta^2}.
  \end{equation}
  Then, for small enough $k$, $\PP\{X_n\to 0\} = 1$.
\end{theorem}

\begin{remark}
  We can treat equation (\ref{eq:nlinIto}) as an equation with noise
  $\xi_n=a+\zeta_n$ where $a=E\xi_n$.
  From this point of view equation (\ref{eq:nlinIto}) is a
  modification of equation (\ref{eq:nonlin_eq}), in which the
  coefficients of the two parts of the noise, drift and diffusion, are
  different. Since $a>0$, the corresponding deterministic equation,
  $$
  x_{n+1} = (1 + kf(x_n)a) x_n+S_n,
  $$
  is unstable.  The diffusion part, $\sqrt{kf(x_n)}\zeta_{n+1}$,
  stabilises the equation.  It becomes possible because the
  coefficient of the diffusion part, $\sqrt{kf(x_n)}$, decreases
  slower then the coefficient of the drift part.

  It is worth noting that equations (\ref{eq:nlinIto}) and
  (\ref{eq:nonlin_eq}) coincide only when $a=0$.
\end{remark}

\section{Proofs}

\subsection{Deterministic lemma}
\label{sec:det_lemma}

For the purposes of comparison with equation (\ref{eq:nlin_recng}), we
discuss here a stability result for the deterministic  equation
\begin{displaymath}
  x_{n+1}=x_n(1+f(x_n)a_{n})+S_n, \quad x_0>0, \quad n\in\N_0.
\end{displaymath}
\begin{lemma}
  \label{lem:deterministic}
  Let $S_n\ge 0$, $f:\R^1 \to [0, 1]$, $f(0)=0$
  and $\inf_{u>c}uf(u)>0$ $\forall c>0$.
  Let also $0>a_n>-1$ and $\sum_{n=1}^{\infty}a_n=-\infty$.

  If $\lim_{n\to \infty} S_n / a_n = 0,$ then $\lim_{n\to
    \infty}x_n=0.$
\end{lemma}
\begin{proof}
  We note that the solution $x_n$ remains positive for all $n$.
  We consider two possibilities: $\liminf x_n>0$ and $\liminf x_n=0$.

  In the first case there exist $c>0$ and $N$ such that $x_n>c$ for all
  $n>N$.  Let $c_1=\inf_{u>c}\{f(u)u\}$ and  $N_1>N$ be such that
  $S_n\le \frac {c_1|a_n|}{2}$ for $n>N_1$.  We have for $n>N_1$
  \begin{displaymath}
    x_{n+1}
    = x_{N_1}+\sum_{i=N_1}^n \left[x_if(x_i)a_{i} + S_i \right]
    \le x_{N_1}
    + c_1\sum_{i=N_1}^n \left[a_{i} + \frac{|a_i|}{2} \right]
    \le x_{N_1} - c_1\sum_{i=N_1}^n \frac{|a_i|}{2}.
  \end{displaymath}
  When $n\to \i$ the right-hand-side of the inequality tends to $-\i$,
  which contradicts the positivity of the solution.  Thus $\liminf x_n=0$.

  Now assume that, even though $\liminf x_n=0$, the lemma is still
  incorrect, i.e.~$\limsup_{n\to \infty}x_n=c>0$.  We fix some $\e<c/2$
  and define
  \begin{displaymath}
    0<\e_1=\inf_{\e<u<2\e}\{f(u)u\}.
  \end{displaymath}
  Now find $N$ such that $S_n < \e_1|a_n|/2$ and $S_n<\e$ whenever
  $n\ge N$.

  If $x_n<\e$ (which must happen infinitely often) with $n>N$,
  we can estimate
  \begin{displaymath}
    x_{n+1}\le x_n+S_n \le 2\e.
  \end{displaymath}
  If, on the other hand, $\e < x_n < 2\e$ then, by
  definition of $\e_1$, $x_nf(x_n) \ge \e_1$ and therefore
  \begin{displaymath}
    x_{n+1} = x_n(1+f(x_n)a_n) + S_n
    \le x_n - \e_1|a_{n}|+\frac {\e_1}{2}|a_{n}|
    < x_{n}.
  \end{displaymath}
  Combining the above two facts we deduce that, once $x_n$ gets
  below $\e$, it cannot increase past $2\e$.  Thus
  $\limsup_{n\to\infty}x_n \le 2\e < c$ and we arrive to a
  contradiction.
\end{proof}

\subsection{Proof of Theorem~\ref{thm:indlin}}
We split the proof into two parts: $\alpha\in(0,1]$ and
$\alpha>1$.

\subsubsection{Proof of Theorem~\ref{thm:indlin} with $\alpha\in(0,1]$}
\label{sec:proof_ii}

We note that $\rho_{i+1}$, defined by
\begin{equation}\label{eq:rho1}
  \rho_{i+1}=X_i^\a(1+\xi_{i+1})^\a-X_i^\a\E(1+\xi_{i+1})^\a
\end{equation}
is an $\F_{n+1}$-martingale-difference.

We apply H\"{o}lder inequality $(x+y)^\alpha \leq x^\alpha +
y^\alpha$ to equation (\ref{eq:lin_recng}) and get
\begin{displaymath}
  \begin{split}
    X_{n+1}^\a &\le  X_n^\a(1+\xi_{n+1})^\a+S_n^\a\\
    &= X_n^\a + X_n^\a\bigl(\E(1+\xi_{n+1})^\a-1\bigr)
    + \bigl[X_i^\a(1+\xi_{n+1})^\a-X_i^\a\E(1+\xi_{n+1})^\a \bigr] + S_n^\a\\
    &= X_n^\a + X_n^\a\bigl(\E(1+\xi_{n+1})^\a-1\bigr) + S_n^\a + \rho_{n+1}\\
    &\le X_n^\a + X_n^\a\bigl[\E(1+\xi_{n+1})^\a-1\bigr]^+ + S_n^\a + \rho_{n+1}
  \end{split}
\end{displaymath}
with $\rho_{n+1}$ defined in equation (\ref{eq:rho1}). We let
\begin{displaymath}
  Y_n = e^{-\sum_{i=1}^n\eta_i}X_n^\a,
  \quad \mbox{with} \quad
  \eta_i=\bigl[\E(1+\xi_{i+1})^\a-1\bigr]^+,
\end{displaymath}
and using the above, arrive at
\begin{displaymath}
  \begin{split}
    Y_{n+1} - Y_n
    &= e^{-\sum_{i=1}^{n+1}\eta_i}\bigl(X_{n+1}^\a-X_n^\a\bigr) +
    X_n^\a e^{-\sum_{i=1}^{n+1}\eta_i}\bigl(1-e^{\eta_{n+1}}\bigr)\\
    &\le e^{-\sum_{i=1}^{n+1}\eta_i}
    \bigl(X_n^\a\eta_{n+1}+\rho_{n+1}+S_n^\a\bigr)
    - \eta_{n+1}X_n^\a e^{-\sum_{i=1}^{n+1}\eta_i}\\
    &= e^{-\sum_{i=1}^{n+1}\eta_i}\rho_{n+1}
    + e^{-\sum_{i=1}^{n+1}\eta_i}S_n^\a
    = \bar\rho_{n+1} + \bar S_n^\a.
  \end{split}
\end{displaymath}
Since $\bar\rho_{n+1}$ is an $\F_{n+1}$-martingale-difference and
$\sum_{i=1}^{\i}\bar S_n^\a<\i$ by a combination of
conditions (\ref{eq:pos_part_conv}) and (\ref{eq:S_alpha_sum}), we can apply
Lemma \ref{lem:nonegdif}. Therefore
$Y_n=\exp\{-\sum_{i=1}^n\eta_i\}X_n^\a$ converges as $n\to \i$.
From condition (\ref{eq:pos_part_conv}) we infer that $X_n^\a$ also \as
converges to a finite limit.

To prove that $\lim_{n\to\infty} X_{n}=0$ we apply
Lemma~\ref{lem:nonegdif} to the inequality
\begin{displaymath}
  X_{n+1}^\a \le
  X_n^\a + X_n^\a\bigl[\E(1+\xi_{n+1})^\a-1\bigr]^-
  + X_n^\a\bigl[\E(1+\xi_{n+1})^\a-1\bigr]^++S_n^\a+\rho_{n+1},
\end{displaymath}
where
$$
\sum_{i=0}^{\i}\bigl[\E(1+\xi_{i+1})^\a-1\bigr]^+X_i^\a
$$
converges \as due to condition (\ref{eq:pos_part_conv}) and the convergence
of $X_n^\a$. From Lemma~\ref{lem:nonegdif} we infer that
$$
-\sum_{i=0}^{\i}\bigl[\E(1+\xi_{n+1})^\a-1\bigr]^-X_i^\a
$$
has to be \as finite.  Combining it with condition (\ref{eq:neg_part_div})
we conclude that $X_i^\a\to0$ \as

\subsubsection{Proof of Theorem~\ref{thm:indlin} with $\alpha>1$}

Let
\begin{equation}
  \label{eq:en}
  \e_n=\frac{\left|1-\E(1+\xi_{n+1})^{\a}\right|}{2\E(1+\xi_{n+1})^{\a}}.
\end{equation}

Applying Lemma \ref{lem:power_ineq} we get
\begin{equation}
  \label{eq:lin_recq}
  X_{n+1}^\a \le (1+\e_n)X_n^\a \Big(1+\xi_{n+1}\Big)^\a +K(\e_n)S_n^\a,
\end{equation}
where $K(\e_n)$  can be estimated by the following
\begin{displaymath}
  \begin{split}
    K(\e_n) &\le 1+K(\a)\e_n^{1-\a}
    = 1 + K(\a)\frac{(2\E(1+\xi_{n+1})^\a)^{\a-1}}
    {\left|1-\E(1+\xi_{n+1})^\a\right|^{\a-1}}\\
    &\le 1 + K(\a)\frac{C^{\a-1}}{\left|1-\E(1+\xi_{n+1})^\a\right|^{\a-1}}
    \le \frac{K_1(\a)}{\left|1-\E(1+\xi_{n+1})^\a\right|^{\a-1}},
  \end{split}
\end{displaymath}
where we used that $\E(1+\xi_{n+1})^\a$ is bounded due to condition
(\ref{eq:pos_part_conv}). From equations (\ref{eq:en})-(\ref{eq:lin_recq}) we
get
\begin{eqnarray*}
  X_{n+1}^\a
  \le X_n^\a &+& X_n^\a
  \left[(1+\e_n)\E\Big(1+\xi_{n+1}\Big)^\a-1\right]\\
  &+& X_n^\a(1+\e_n)\left[\Big(1+\xi_{n+1}\Big)^\a
  - \E\Big(1+\xi_{n+1}\Big)^\a\right] + K(\e_n)S_n^\a.
\end{eqnarray*}
By substituting the value of $\e_n$ into
$(1+\e_n)\E\left(1+\xi_{n+1}\right)^\a-1$ we see that it is equal
to $\left[\E\left(1+\xi_{n+1}\right)^\a-1\right]/2$ when
$\E\left(1+\xi_{n+1}\right)^\a<1$ and to
$3\left[\E\left(1+\xi_{n+1}\right)^\a-1\right]/2$ otherwise.  That
is, we can write
\begin{displaymath}
  (1+\e_n)\E\Big(1+\xi_{n+1}\Big)^\a-1
  = \frac12\left[\E\Big(1+\xi_{n+1}\Big)^\a-1\right]^-
  + \frac32\left[\E\Big(1+\xi_{n+1}\Big)^\a-1\right]^+.
\end{displaymath}
We finally arrive to

\begin{eqnarray*}
  X_{n+1}^\a \le X_n^\a
  &+& \frac12 X_n^\a \left[\E\Big(1+\xi_{n+1}\Big)^\a-1\right]^-
  + \frac32 X_n^\a \left[\E\Big(1+\xi_{n+1}\Big)^\a-1\right]^+ \\
  &+& \rho_{n+1}
  + \frac{K_1(\a)}{(1-\E(1+\xi_{n+1})^{\a})^{\a-1}}S_n^\a,
\end{eqnarray*}
where $\rho_{n+1}$ is an $\F_{n+1}$-martingale-difference. Now we
apply Lemma \ref{lem:nonegdif} and complete the proof as in
Section~\ref{sec:proof_ii}.

\subsection{Proof of Theorem~\ref{thm:indexp}}

We mimic the proof of Theorem~\ref{thm:indlin} (see
Section~\ref{sec:proof_ii}) with
$Y_n=e^{-\sum_{i=1}^n\kappa_i}X_n^\a$ to get
\begin{displaymath}
  Y_{n+1}-Y_n \le \bar\rho_{n+1}+ e^{-\sum_{i=1}^{n+1}\kappa_i}S_n^\a.
\end{displaymath}
Because $\bar\rho_{n+1}$ is an $\F_{n+1}$-martingale-difference
and due to condition~(\ref{eq:kappa1}), we can apply Lemma
\ref{lem:nonegdif}.  Hence we get that
$Y_n=\exp\{-\sum_{i=1}^n\kappa_i\}X_n^\a$ converges to a finite limit as
$n\to \i$.  Then for every $\gamma\in(0, 1)$
\begin{equation}
  \exp\left\{-\gamma\sum_{i=1}^n\kappa_i\right\}X_n^\a
  \le Y_n \exp\left\{(1-\gamma)\sum_{i=1}^n\kappa_i\right\} \to 0
\end{equation}
using condition~(\ref{eq:kappa2}).

\subsection{Proof of Lemma~\ref{lem:limsup_counter}}

We choose $\gamma$ such that $\alpha<\gamma<\beta$ and take $S_n =
n^{-1/\gamma}$ so that condition (\ref{eq:limsup_betaSum}) is clearly
satisfied.  Now define the distribution of $\xi_n$ so that $1+\xi_n$
takes values in $(a,\infty)$, $a>0$, with the density function
\begin{displaymath}
  p(x) = \frac{\gamma a^\gamma}{x^{1+\gamma}}.
\end{displaymath}
First we ascertain that $\E(1+\xi)^\alpha = 1$.  Indeed, since
$\alpha<\gamma$,
\begin{equation*}
  \E(1+\xi)^\alpha
  = \gamma a^\gamma \int_a^\infty x^{-1-\gamma+\alpha} dx
  = \frac{\gamma}{\gamma-\alpha}a^\alpha
\end{equation*}
and condition~(\ref{eq:limsup_Ealpha}) can be satisfied with an appropriate
choice of $a$.

Now we can study the behaviour of solutions of equation~(\ref{eq:lin_recng}).
Since both summands in the right hand side of equation~(\ref{eq:lin_recng})
are positive, $X_{n+1} \geq S_n$ and therefore $X_{n+2} \geq
(1+\xi_{n+2})S_n$.  Define the sequence of independent events $A_n =
\left\{(1+\xi_{n+2})S_n > C\right\}$, where $C>0$ is an arbitrary constant.
We have
\begin{displaymath}
  \PP(A_n) = \PP\left((1+\xi_{n+2}) > \frac{C}{S_n}\right)
  = \PP\left((1+\xi_{n+2}) > Cn^{1/\gamma}\right)
  = \frac{a^\gamma}{C^\gamma} n^{-1}.
\end{displaymath}
Thus,
\begin{displaymath}
  \sum_{n=1}^\infty \PP(A_n) = \infty,
\end{displaymath}
and, by Borel-Cantelli lemma, events $A_n$ must happen infinitely often.
Therefore, infinitely often $X_n > C$.  Since $C$ was arbitrary, we conclude
that $\limsup_{n\to\infty} X_n = \infty$ \as

\subsection{Proof of Theorem~\ref{thm:liminf}}

Assume the contrary, for some $s>0$ the event $J_s = \{\omega \colon
\inf_{n} X_n^\alpha > s \}$ occurs with non-zero probability.  Fix
$\epsilon>0$ such that $\alpha\E\ln(1+\xi)+\ln(1+\epsilon)<0$ and
consider the event $\Theta =\{\omega \colon \sum_{i=1}^n
\ln\big((1+\epsilon)(1+\xi_i)^\alpha\big) \to -\infty \}$.  By
applying the law of large numbers it is straightforward to show that
$\Theta$ occurs with probability 1.

Raising recursion (\ref{eq:lin_recng}) to power $\alpha$ we get by
Lemma~\ref{lem:power_ineq}
\begin{equation}
  \label{eq:lin_rec_alpha}
  X_{n+1}^\a \le (1+\epsilon)X_n^\a \Big(1+\xi_{n+1}\Big)^\a 
  + K(\epsilon)S_n^\a.
\end{equation}
Now let $n$ be such that $K(\epsilon) S_n^\alpha < s/2$.  Restricting
our attention to $\omega \in J_s$ we apply logarithm to both sides of
inequality (\ref{eq:lin_rec_alpha}) and use the inequality
\begin{displaymath}
  \ln(x+y) \leq \ln(x) + \frac{y}{x}
\end{displaymath}
to obtain
\begin{displaymath}
  \ln X_{n+1}^\alpha
  \leq \ln\Big(X_n^\alpha (1+\epsilon)(1+\xi_{n+1})^\alpha\Big)
  + \frac{K(\e)S_n^\alpha}{X_n^\alpha(1+\epsilon)(1+\xi_{n+1})^\alpha}.
\end{displaymath}
Combining inequality (\ref{eq:lin_rec_alpha}) and the definition of $J_s$ we
can estimate $X_n^\alpha (1+\epsilon) (1+\xi_{n+1})^\alpha \geq X_{n+1}^\alpha
- K(\epsilon) S_n^\alpha > s/2$ and, therefore,
\begin{displaymath}
  \ln X_{n+1}^\alpha \leq \ln(X_n^\alpha)
    + \ln\Big((1+\epsilon) (1+\xi_{n+1})^\alpha\Big)
    + K(\epsilon)\frac{S_n^\alpha}{s/2}.
\end{displaymath}

Applying the above inequality recursively we obtain
\begin{displaymath}
  \ln X_{n+k}^\alpha <  \ln(X_n^\alpha)
  + \sum_{i=1}^k \ln\Big((1+\epsilon) (1+\xi_{n+i})^\alpha\Big)
  + C\sum_{i=0}^{k-1} S_{n+i}^\alpha,
\end{displaymath}
where $C=2K(\epsilon)/s$.  Since $X_{n+k}^\alpha>s$ and $S_n^\alpha$ are
summable to, say, $S$, we conclude that for all $k$
\begin{displaymath}
  \sum_{i=1}^k \ln\Big((1+\epsilon) (1+\xi_{n+i})^\alpha\Big)
  > \ln(s) - \ln(X_n^\alpha) - CS
\end{displaymath}
and, therefore, the event $\omega$ cannot belong to $\Theta$. Thus
$J_s\cap \Theta = \emptyset$ which is a contradiction.

\subsection{Proof of Lemma~\ref{lem:alnind}}

Taking the expectation of the Taylor expansion of $(1+\xi_i)^\a$
in terms of $\alpha$ we get
\begin{displaymath}
  \E(1+\xi_{i})^{\a} = 1 + \a\E\ln(1+\xi_{i})
  + \a^2\E\left(\frac{\ln^2(1+\xi_{i})}{2}(1+\xi_{i})^{\theta}\right),
\end{displaymath}
where $\theta\in[0,\a]$.  The left side of estimate (\ref{eq:uniform_bound})
is then obtained by leaving out the third term.

To estimate $\E\left(\ln^2(1+\xi_{i})(1+\xi_{i})^\theta/2\right)$ from
above we consider two cases: $1+\xi_{i}>1$ and $1+\xi_{i}<1$.
Since $\theta\leq\a\leq1$, in the first case we have
$(1+\xi_{i})^{\theta}\le (1+\xi_{i})$, while in the second
$(1+\xi_{i})^{\theta}\le (1+\xi_{i})^0=1$. Then, in both cases, we have
\[
 (1+\xi_{i})^{\theta}\le 2+\xi_i.
\]
If $\E\ln(1+\xi_{i})$ is negative we continue with
\begin{eqnarray*}
  \E(1+\xi_{i})^{\a}
  &\le& 1 + \a\E\ln(1+\xi_{i})\biggl(1-\a \frac{\E \left((2+\xi_{i})\ln^2
      (1+\xi_{i})\right)}{2\left|\E \ln(1+\xi_{i})\right|}\biggr)\\
  &\le& 1 + \a\E\ln (1+\xi_{i})\biggl(1-\a\frac{K}{2}\biggr)
  \le 1+\frac{\a}{2} \E\ln(1+\xi_{i}),
\end{eqnarray*}
while if $\E\ln(1+\xi_{i})>0$ we obtain by a similar calculation
\[
\E(1+\xi_{i})^{\a}\le 1+\a \frac{3\E\ln (1+\xi_{i})}{2}.
\]

\subsection{Proof of Theorem~\ref{thm:posnegn}}

We note that $\rho_{i+1}$, defined by
$$
\rho_{i+1}=f(X_i)X_i\xi_{i+1}-f(X_i)X_i\E(\xi_{i+1}).
$$
is an $\F_{n+1}$-martingale-difference.

After rearranging in equation (\ref{eq:nonlin_eq}) we get recursively
\begin{equation}
  \begin{split}
    \label{eq:nEyn2}
    X_{n+1} &= X_n + f(X_n)X_n\E\xi_{n+1}
    + \bigl[ f(X_n)X_n\xi_{n+1}- f(X_n)X_n\E\xi_{n+1}\bigr]+S_n\\
    &= X_n + f(X_n)X_n\bigl[\E\xi_{n+1}\bigr]^+
    + f(X_n)X_n\bigl[\E\xi_{n+1}\bigr]^- + \rho_{n+1} + S_n\\
    &\le  X_n + f(X_n)X_n\bigl[\E\xi_{n+1}\bigr]^+ + \rho_{n+1} + S_n\\
    &\le  X_n + X_n\bigl[\E\xi_{n+1}\bigr]^+ + \rho_{n+1} + S_n.
  \end{split}
\end{equation}
From this point we continue as in Section~\ref{sec:proof_ii} with
$\eta_i=[\E\xi_{n+1}\bigr]^+$ and conclude that $X_i$ converges to
a finite limit a.s.  Then
\begin{displaymath}
  \sum_{i=0}^{\i}[\E(\xi_{i+1})]^+f(X_i)X_i
\end{displaymath}
is \as finite.  Applying Lemma \ref{lem:nonegdif} again (to the
second line in inequality (\ref{eq:nEyn2})), we conclude that
\begin{displaymath}
  \sum_{i=0}^{\i}[\E(\xi_{i+1})]^-f(X_i)X_i
\end{displaymath}
also has to be \as finite.  If condition (\ref{eq:neg_exp_div}) is
fulfilled, $f(X_i)X_i$ is forced to converge to zero.  Therefore
$X_i\to0$ \as

\subsection{Proof of Theorem~\ref{thm:asumind}}

Applying the inequality
\begin{equation}
  \label{eq:ser_alpha}
  (1+x)^\alpha \le 1 + \alpha x - \frac{1-\alpha}{2} x^2
  + \frac{(1-\alpha)(2-\alpha)}{6} x^3, \qquad x>-1,\ 0<\alpha<1
\end{equation}
and noting that $f^2(X_n)\ge f^3(X_n)$, we obtain
\begin{equation*}
  \begin{split}
    &\E\biggl[X_n^\a(1+f(X_n)\xi_{n+1})^\a\bigl| \F_{n}\biggr]\\
    &\le X_n^\a \left(1 + \a f(X_n)\E\xi_{n+1}
    - \frac{1-\a}{2}f^2(X_n)\E\xi_{n+1}^2
    + \frac{(1-\a)(2-\a)}{6}f^3(X_n)\E\xi_{n+1}^3 \right)\\
    &\le X_n^\a + \a X_n^\a f(X_n)[\E\xi_{n+1}]^+
    - \frac{1-\a}{2}X_n^\a f^2(X_n)
    \left(\E\xi_{n+1}^2 - \frac{(2-\a)}{3}[\E\xi_{n+1}^3]^+\right).
  \end{split}
\end{equation*}
Now, applying inequality $(a+b)^\a\le a^\a+b^\a$
to equation~(\ref{eq:nonlin_eq}) we get from the above
\begin{displaymath}
  \begin{split}
    X_{n+1}^\a &\le  X_n^\a \Big(1+f(X_n)\xi_{n+1}\Big)^\a + S_n^\a\\
    &= \E\biggl[X_n^\a(1+f(X_n)\xi_{n+1})^\a\bigl| \F_{n}\biggr]\\
    &\hspace{1cm} + \biggl(X_n^\a(1+f(X_n)\xi_{n+1})^\a
    - \E\biggl[X_n^\a(1+f(X_n)\xi_{n+1})^\a
    \bigl|\F_{n}\biggr]\biggr)+ S_n^\a\\
    &\le X_n^\a + \a X_n^\a f(X_n)[\E\xi_{n+1}]^+ \\
    & \hspace{1cm} - \frac{1-\a}{2}X_n^\a f^2(X_n)\left(\E\xi_{n+1}^2
      -\frac{(2-\a)}{3}[\E\xi_{n+1}^3]^+\right) + \rho_{n+1}+S_n^\a.
  \end{split}
\end{displaymath}
Now we complete the proof in the same way as in
Theorem~\ref{thm:posnegn}.

\subsection{Proof of Theorem~\ref{thm:div}}

For the proof we need some preliminary facts.
\begin{lemma}
  \label{lem:mart}
  Let $\{X_n\}_{n\in \N}$ be a sequence of $\F_n$-measurable
  random variables such that ${\bf E}\bigl(X_n\bigl| {\cal
    F}_{n-1}\bigr)= 1$. Let $Z_n=\prod_{i=1}^nX_i$ and $\E |Z_n|<\i$
  for all $n\in\N$.  Then $\{Z_n\}_{n\in \N}$ is a martingale.
\end{lemma}
\begin{proof}
  To check the martingale condition ${\bf E}\bigl(Z_n\bigl| {\cal
    F}_{n-1}\bigr)= Z_{n-1}$ we use the $\F_{n-1}$-measurability
  of $Z_{n-1}$:
  \begin{displaymath}
    \E\bigl(Z_n\bigl| \F_{n-1}\bigr)
    = \E\bigl(Z_{n-1}X_n\bigl| \F_{n-1}\bigr)
    = Z_{n-1}\E\bigl(X_n\bigl| \F_{n-1}\bigr) = Z_{n-1}.
  \end{displaymath}
\end{proof}

\begin{lemma}\label{lem:martrepr}
  Let $X_n$ be a solution of equation~(\ref{eq:nonlin_eq}). Then the sequence
  $\{M_n\}_{n\in \N}$, defined by
  \begin{equation}
    \label{eq:Mn}
    M_n = \prod_{i=0}^{n-1}\frac{(1 +f(X_i)\xi_{i+1})^{-1}}
    {\E\left((1 + f(X_i)\xi_{i+1})^{-1}\bigl|\F_{i}\right)}
  \end{equation}
  is an $\F_n$- martingale.
\end{lemma}
\begin{proof}
  To make sure that our definition makes sense we estimate
  \begin{equation}
    \label{eq:twoside}
    1 + f(X_i)\xi_{i+1} \ge 1-|\xi_{i+1}| > 1-k_0 > 0,
  \end{equation}
  therefore $\E\left((1 +
  f(X_i)\xi_{i+1})^{-1}\bigl|\F_{i}\right)$ is well defined.
  Because $M_n$ is always positive, we can write $\E|M_n| = \E M_n
  = \E M_1 = 1 < \infty$.
  Now we apply Lemma~\ref{lem:mart} to conclude the proof.
\end{proof}

The lemma below is a variant of the theorem of convergence of
non-negative martingale (see e.g. \cite{LSh}).
\begin{lemma}
  \label{lem:posmart}
  If $\{X_n\}_{n\in \N}$ is non-negative martingale, then
  $\lim_{n\to \i} X_n$ exists with probability 1.
\end{lemma}
From Lemma \ref{lem:martrepr} and Lemma \ref{lem:posmart} we can
get
\begin{corollary}
  \label{lem:martconv}
  Let $\{M_n\}_{n\in \N}$ be the martingale defined by (\ref{eq:Mn}), then
  $\lim_{n\to \i} M_n$ exists with probability 1.
\end{corollary}

Now we proceed to the proof of the theorem.  First we note that the
solution $X_n$ of equation~(\ref{eq:nonlin_eq}) can be represented in the
following form
\begin{equation}\label{eq:reprx}
  X_{n} = X_0 M_n^{-1}\prod_{i=0}^{n-1}
  \frac1{\E\left((1 + f(X_i)\xi_{i+1})^{-1}\bigl| \F_{i}\right)}.
\end{equation}
Here $M_n$ is defined by equation (\ref{eq:Mn}) and, by Corollary
\ref{lem:martconv}, $M_n\le H_1$ with {a.s.} finite random
variable $H_1=H_1(\omega)$.

Suppose now that theorem is not correct.  Then there exists a set
$\Omega_1\subseteq \Omega$ of non-zero probability such that $X_n\to 0$ \as on
$\Omega_1$. We aim to show that for any $\omega\in\Omega_1$, there exists
$N(\omega)$ such that
\begin{displaymath}
  \E\biggl(\big(1 + f(X_i)\xi_{i+1}\big)^{-1}\bigl| {\cal
  F}_{i}\biggr)\le 1, \qquad \forall i\ge N(\omega).
\end{displaymath}
For $\forall i\in \N$ we can perform the Taylor expansion
\begin{displaymath}
  (1 + f(X_i)\xi_{i+1})^{-1}
  = 1 - f(X_i)\xi_{i+1} + f^2(X_i)\xi^2_{i+1}
  -\frac{f^3(X_i)\xi^3_{i+1}}{(1+\theta_{i+1})^4}
\end{displaymath}
with $\theta_{i+1}$ lying between 0 and $f(X_i)\xi_{i+1}$. Using
equation (\ref{eq:twoside}) and noting that $X_n$ is positive and
\begin{displaymath}
  \E\bigl(f(X_i)\bigl| \F_i\bigr) = f(X_i),
  \quad
  \E\bigl(\xi_{i+1}\bigl| \F_i\bigr) =\E \xi_{i+1},
  \quad 0 \le f(X_i) \le 1,
\end{displaymath}
we estimate
\begin{displaymath}
  \E\left(\frac{f^3(X_i)\xi^3_{i+1}}{(1+\theta_{i+1})^4}\biggl| \F_{i}\right)
  \le \frac{L^3 f^3(X_i)}{(1-k_0)^4}.
\end{displaymath}
Then we have
\begin{eqnarray*}
  \E\biggl(\big(1 + f(X_i)\xi_{i+1}\big)^{-\a}\biggl| \F_i\biggr)
  &\le& 1 - f(X_i)\E\xi_{i+1} + f^2(X_i)L^2
  + \frac{L^3f^3(X_i)}{(1-k_0)^4}\\
  &=& 1 - f(X_i)\left(\E\xi_{i+1} - f(X_i)L^2
    - \frac{L^3f^2(X_i)}{(1-k_0)^4}\right).
\end{eqnarray*}
The function $f$ is such that $f(X_n)\to0$ \as on $\Omega_1$, therefore we
can find such $N(\omega)$ that for $\Omega_1$ and $i\ge N(\omega)$
\begin{displaymath}
  \E\biggl((1 + f(X_i)\xi_{i+1})^{-\a}\biggl| \F_i\biggr)
  \le 1 - f(X_i)\frac{\E\xi_{i+1}}{2} < 1.
\end{displaymath}
Combining this with representation (\ref{eq:reprx}) and with \as boundedness
of $M_n$ we conclude that solution $X_n$ cannot tend to 0 on $\Omega_1$.

\subsection{Proof of Theorem~\ref{thm:Ito}}

As before, we raise equation (\ref{eq:nlinIto}) to power $\alpha$ and set
\begin{displaymath}
  \rho_{n+1}
  = X_n^\alpha \left(1 + k f(X_n)a + \sqrt{kf(X_n)}\zeta_{n+1}\right)^{\a}
   - \E\left[\left. X_n^\alpha\left(1 + k f(X_n)a +
        \sqrt{kf(X_n)}\zeta_{n+1}\right)^{\a}
    \right| \F_n \right].
\end{displaymath}
We now aim to show that the conditional expectation above is
negative. Applying inequality~(\ref{eq:ser_alpha}) and remembering
that $\E\zeta_{n+1}=0$, we get
\begin{displaymath}
  \begin{split}
    &\E\left[\left. \left(1 + k f(X_n)a +
          \sqrt{kf(X_n)}\zeta_{n+1}\right)^{\a}
      \right| \F_n \right] \\
    &\hspace{1cm} \le 1 + \alpha kf(X_n) a
    - \frac{\alpha(1-\alpha)}2 \left( (akf(X_n))^2 + kf(X_n)\E\zeta^2
    \right)\\
    &\hspace{2cm} + \frac{\alpha(1-\alpha)(2-\alpha)}6 \left( (akf(X_n))^3
      + 3a(kf(X_n))^2\E\zeta^2 + (kf(X_n))^{3/2}\E\zeta^3 \right)\\
    &\hspace{1cm} \le 1 + \alpha kf(X_n) \left( a
      - \frac{1-\alpha}2\E\zeta^2 + O(\sqrt{k}) \right)
  \end{split}
\end{displaymath}
Due to condition (\ref{eq:Ito_alpha}) there exist $k_0$ and $a_0$,
such that for $k<k_0$
\begin{displaymath}
  a - \frac{1-\a}2\E\zeta^2 + O(\sqrt{k}) \leq -a_0<0.
\end{displaymath}
Therefore we obtain the estimation
\begin{displaymath}
  X_{n+1}^{\a} \le X_{n}^{\a} - a_0 X_{n}^{\a}\a k f(X_n) +
  \rho_{n+1} + S_n^{\a}.
\end{displaymath}
Now we can apply Lemma~\ref{lem:nonegdif} and complete the proof
by the familiar method.

\section*{Acknowledgement}
The authors are grateful to D.Cline for his helpful comments,
to J.Appleby for the interesting discussion of the similar problem
in the continuous case, to B.Winn for critically reading the
manuscript and pointing out various misprints and to a referee for
useful suggestions.


\end{document}